\documentclass[11pt]{article}

\usepackage[latin1]{inputenc}
\usepackage{subfigure}
\usepackage{latexsym}
\usepackage{amsfonts}

\usepackage{epsfig}
\usepackage{amsmath,amssymb,stmaryrd,eufrak,float}
\usepackage{graphics, graphicx}
\usepackage{color}
\usepackage[all]{xy}

\newtheorem{proposition}{Proposition}

\author{E. Grazzini$^{1}$ \and E. Pergola$^{1}$
\and M. Poneti$^{2}$}

\title{On the exhaustive generation of convex permutominoes}

\date{}
\begin{document}

\hyphenation{critical permu-tomino Moreover parallelogram
obtained}
 \maketitle
 \footnotetext[1]{Dipartimento di Sistemi e Informatica, Viale G. B.
Morgagni 65, 50134 Firenze, Italy,  {\tt
\{ely,elisa\}@dsi.unifi.it }}

\footnotetext[2]{Dipartimento di Scienze Matematiche e
Informatiche R. Magari, Pian dei Mantellini 44, 53100 Siena,
Italy, {\tt poneti@unisi.it }\\}
%\begin{description}
%\item[$^a$] Dipartimento di Sistemi e Informatica, Viale G. B.
%Morgagni 65, 50134 Firenze, Italy,  {\tt \{ely,elisa\}@dsi.unifi.it }\\
 %\item[$^b$]Dipartimento di Scienze Matematiche e Informatiche R. Magari, Pian dei Mantellini 44, 53100 Siena, Italy,
%{\tt poneti@unisi.it }\\

 %\end{description}

\begin{abstract}
A permutomino of size $n$ is a polyomino determined by a pair
$(\pi _1 ,\pi _2)$ of permutations of size $n+1$, such that
$\pi_1(i) \neq \pi_2(i)$, for $ 1 \leq i \leq n+1$. In this paper,
after recalling some enumerative results about permutominoes, we
give a first algorithm for the exhaustive generation of a
particular class of permutominoes, the {\em convex} permutominoes,
proving that its cost is proportional to the number of generated
objects.
\end{abstract}

\section{Introduction}
A permutomino is a special polyomino, defined by two permutation matrices having
the same size. The class of permutominoes  was introduced by Incitti in
\cite{incitti} while studying the problem of determining the
$\widetilde{R}$-polynomials (related to the Kazhdan-Lusztig
R-polynomials) associated with a pair $(\mu,\nu)$ of
permutations. In his paper Incitti gave a general definition of
these combinatorial objects using some algebraic notions.

%A new definition of permutominoes was first given in \cite{mathinfo}.
In this paper we use the definition of permutominoes given in \cite{dfpr} which does not use the algebraic
notions and nevertheless, though different, it turns out to be
equivalent to Incitti's one.

The main results about permutominoes concern the enumeration of
various subclasses of permutominoes and the characterization for
the permutations defining these subclasses
\cite{bdpr,boldi,dfpr,mathinfo,new},
%but little is concerned with
while at our knowledge nothing exists about their generation. On
the other hand, exhaustive generation of combinatorial objects
\cite{acta04,puma,acta07} is an area of increasing interest. In
fact, many practical questions in diverse areas, such as hardware
and software testing, and combinatorial chemistry, require for
their solution the exhaustive search through all objects in the
class.

Actually, in \cite{dfpr} a recursive generation of all convex
permutominoes of size $(n+1)$ from the ones of size $n$, according
to the ECO method \cite{eco}, is presented. Section \ref{basic}
contains basic definitions and some enumerative results of convex
permutominoes of size $n$. Section \ref{eco} recalls the recursive
generation of convex permutominoes presented in \cite{dfpr} and
Section \ref{algo} illustrates the exhaustive generating algorithm
based on the recursive construction recalled in Section \ref{eco}.

\section{Basics on permutominoes}\label{basic}
\subsection {Definitions and properties}

In order to define permutominoes we need to introduce polyominoes.
In the plane $\Bbb Z \times \Bbb Z$ a {\em cell} is a unit square
and a {\em polyomino} is a finite connected union of cells having
no cut point. Polyominoes are defined up to translations. A {\em
column} ({\em row}) of a polyomino is the intersection between the
polyomino and an infinite strip of cells whose centers lie on a
vertical (horizontal) line.

\begin{figure}[htb]
\centerline{\hbox{
\includegraphics[width=0.7\textwidth]{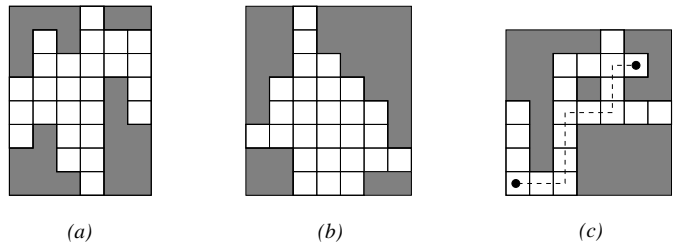}
}} \caption{$(a)$ a column convex polyomino; $(b)$ a convex
polyomino; $(c)$ a directed (not convex) polyomino with a
hole.\label{dirconv}}
\end{figure}

In order to simplify many problems which are still open on the
class of polyominoes, several subclasses were defined by combining
two notions: the geometrical notion of {\em convexity} and the
notion of {\em directed growth}.
%, which comes from statistical physics.
A polyomino is said to be {\em column convex} [{\em
row convex}] if its intersection with any vertical [horizontal]
line is convex (Figure~\ref{dirconv}~$(a)$). A polyomino is {\em
convex} if it is both column and row convex
(Figure~\ref{dirconv}~$(b)$). In a convex polyomino the {\em
semi-perimeter} is given by the sum of the number of rows and
columns, while the {\em area} is the number of its cells.

A polyomino $P$ is said to be {\em directed} when every cell of
$P$ can be reached from a distinguished cell, called {\em
root} (usually the bottom leftmost cell), by a path
which is contained in $P$ and uses only north and east unit steps
(Figure~\ref{dirconv}~$(c)$).
%Figure~\ref{par1}~$(d)$ depicts a polyomino which is both directed and convex.

Let $P$ be a polyomino without holes having $n$ rows and $n$
columns, $n\geq 1$; without loss in generality, we assume that the
bottom leftmost vertex of the polyomino minimal
bounding rectangle lies in $(1,1)$.
Let $A_1, \ldots , A_{2(r+1)}$ be the sequence of
the vertices of $P$ obtained by visiting the boundary in clockwise
sense, starting from its leftmost point with minimal
ordinate.

\begin{figure}[h]
\begin{center}
\centerline{\hbox{\includegraphics[width=0.9\textwidth]{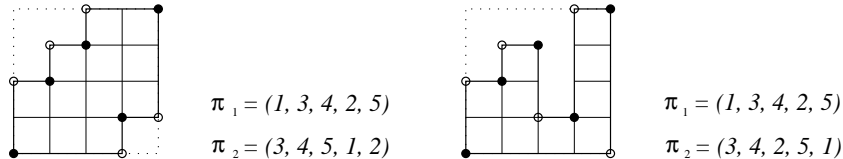}}}
\caption{Two permutominoes and the associated permutations. The
permutation $\pi_1$ (resp. $\pi_2$) is represented by black (resp.
white) dots. }\label{sample}
\end{center}
\end{figure}

We say that $P$ is a {\em permutomino} if the sets ${{\cal P}_1}=
\{A_1, A_3,\ldots ,A_{2r+1}\}$ and ${{\cal P}_2}= \{A_2,
A_4,\ldots ,A_{2r+2}\}$ represent two permutations matrices of
[$n+1$] = \{1,2,\ldots , $n+1$\}. Obviously, if $P$ is  a
permutomino, then $r =n$, and $n$ is called the {\em size} of the
permutomino. The two permutations defined by ${{\cal P}_1}$ and
${{\cal P}_2}$ are indicated by $( \, \pi _1(P) , \, \pi _2
(P)\,)$ (briefly, $(\pi _1 ,\pi _2)$), respectively (see Figure
~\ref{sample}). Given a pair of
permutations $\omega = (\mu,\,\nu)$ of $[n+1]$, we say that a
permutomino $P$ is {\em associated} with  $\omega$  if $\mu =
\pi_1(P)$ and $\nu = \pi_2(P)$.

A permutomino is {\em convex} ({\em directed}) if it is a convex (directed) polyomino.
A {\em parallelogram} permutomino is a directed and convex one having the $(1,1)$ and
$(n,n)$ vertices in common with its minimal bounding square; a {\em stack} permutomino
is a directed and convex one in which the bottom side of its minimal bounding square belongs
to the permutomino itself.

The definition of permutominoes leads to the following remarkable
property:

\begin{proposition}\label{prop1}
Any permutomino $P$ has the property that, for each abscissa (ordinate) there exists
exactly one vertical (horizontal) side in the boundary of $P$ having
such coordinate. This property is also a sufficient condition for a polyomino to be
a permutomino.
\end{proposition}

\begin{figure}[h]
\begin{center}
%\centerline{\hbox{\includegraphics[width=0.9\textwidth]{operationA.eps}}}
\centerline{\hbox{\includegraphics[height=4cm]{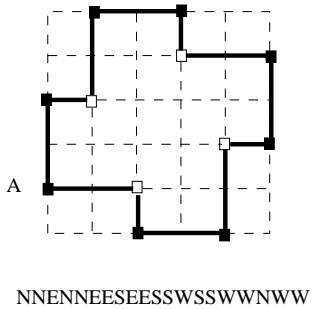}}}
\caption{The coding of the boundary of a permutomino, starting
from A and moving in clockwise sense; its salient (resp.
reentrant) points are indicated by black (resp. white)
squares.}\label{sample2}
\end{center}
\end{figure}

Starting from the leftmost point having minimal ordinate, and
moving in a clockwise sense, the boundary of a permutomino $P$ can
be encoded as a word in a four letter alphabet, $\{N,E,S,W\}$,
where $N$ (resp., $E$, $S$, $W$) represents a {\em north} (resp.
{\em east}, {\em south}, {\em west}) unit step. Any occurence of a
sequence $NE$, $ES$, $SW$ or $WN$ in the word encoding $P$ defines
a {\em salient point} of $P$, while any occurence of a sequence
$EN$, $SE$, $WS$ or $NW$ defines a {\em reentrant point} of $P$
(see Figure~\ref{sample2}). For simplicity of notation and to
clarify the definition of the construction recalled in Section
~\ref{eco}, the reentrant points of a convex permutomino are
grouped in four classes; in practice, the reentrant point
determined by a sequence {\em EN} (resp. {\em SE}, {\em WS}, {\em
NW}) is represented with the symbol $\alpha$ (resp. $\beta$,
$\gamma$, $\delta$).

%******* INSERIRE FIGURA (DA FARE)******

\subsection{Previous enumerative results}

Let us recall the main enumerative results concerning convex
permutominoes. In \cite{mathinfo, new}, using bijective
techniques, the authors provide enumeration of various classes of
convex permutominoes, including the {\em parallelogram}, the {\em
directed convex} and the {\em stack} ones; moreover, a
characterization of the permutations associated with permutominoes
of each class is given. Let ${\cal C}_n$  (resp. ${\cal P}_n$,
${\cal D}_n$, ${\cal S}_n$)    be the set of convex (resp.
parallelogram, directed convex, stack) permutominoes of size $n$
and
$$
\widetilde{\cal C}_n = \left \{ \, \pi_1(P) \, : \, P \in {\cal
C}_n \, \right \}, \quad \widetilde{\cal P}_n = \left \{ \,
\pi_1(P) \, : \, P \in {\cal P}_n \, \right \},
$$
$$
\widetilde{\cal D}_n = \left \{ \, \pi_1(P) \, : \, P \in {\cal
D}_n \, \right \}, \quad \widetilde{\cal S}_n = \left \{ \,
\pi_1(P) \, : \, P \in {\cal S}_n \, \right \}.
$$
The enumeration results obtained in \cite{mathinfo, new} are shown in Table~\ref{ric}.
\begin{table}[h]
\centerline{
\begin{tabular}{c|c|c}
  $Class\mbox{ }of\mbox{ }convex\mbox{ }$ &$Number\mbox{ } of\mbox{ } convex$ &$Associated\mbox{ }set$\\
  $permutominoes$ & $permutominoes\mbox{ }with\mbox{ }size\mbox{ }n$&$of\mbox{ }permutations$\\
  \hline
  \hline
  $$ & & \\
  $Parallelogram$ &$| {\cal P}_n| = c_{n}$ &$| \widetilde {\cal P}_{n}| = c_{n-1}$ \\
  $$ & & \\
  $Directed\mbox{ }convex$ &$| {\cal D}_n| = \frac{1}{2}b_{n}$ & $| \widetilde {\cal D}_{n}| = b_{n-1}$ \\
  $$ & & \\
  $Stack$ &$| {\cal S}_n| = 2^{n-1}$ & $| \widetilde {\cal S}_{n}| =2^{n-1} $ \\
   \end{tabular}}\caption{Enumeration results.}\label{ric}
\end{table}

\noindent where $c_n$ is the $n$th {\em Catalan number},

$$ c_n = \frac{1}{n+1}{2n \choose n}$$
and $b_n$ are the {\em central binomial coefficients}

$$b_n ={2n \choose n}.$$

In \cite{dfpr} it was proved, using the ECO method \cite{eco},
that the number of {\em convex permutominoes} of size $n$ is:

$$ 2(n+3)4^{n-2} -\frac{n}{2}{2n \choose n}, \quad n \geq 1.$$
The first terms of the sequence are
\begin{center}
$1,4,18,84,394,1836,8468, \dots$
\end{center}
(sequence A126020 in \cite{sloane}). The same enumerative result was obtained, using a different
approach, by Boldi et al. \cite{boldi}.\\
An interesting study of combinatorial properties of the set $\widetilde{\cal C}_n $ is in \cite{bdpr}.

\section{ECO construction of convex permutominoes}\label{eco}

In this section we recall the  ECO construction of convex
permutominoes as given in \cite{dfpr}.

Let ${\cal C}_n$ be the set of convex
permutominoes of size $n$ and let $P \in {\cal C}_n$; the number of cells in the rightmost
column of $P$ is called the {\em degree} of $P$.
Let us consider the following properties of a convex permutomino:
\begin{description}
\item{\textbf{U1}} : the uppermost cell of the rightmost column of
$P$ has the maximal ordinate among all the cells of the
permutomino;
 \item{\textbf{U2}}  : the lowest cell of the rightmost
column of $P$ has the minimal ordinate among all the cells of the
permutomino.
\end{description}
According to the ECO method \cite{eco}, it is necessary to define
an operator $\vartheta:~{\cal C}_n~\rightarrow~2^{{\cal C}_{n+1}}$
which defines a recursive construction of all the convex
permutominoes of size $(n+1)$ in a unique way from the objects of
size $n$. The operator $\vartheta$ defined in \cite{dfpr} acts on
a convex permutomino performing some local expansions on the cells
of its rightmost column. Let $c_1, \ldots, c_n$ (resp. $r_1,
\ldots, r_n$) be the columns (resp. rows) of a permutomino $P$ of
size $n$ numbered from left to right (resp. bottom to top) and let
$\ell(c_i)$ (resp. $\ell(r_i)$) be the number of cells in the
$i$th column (resp. row), with $1 \leq i \leq n$. The four
operations of $\vartheta$, denoted by $(\alpha), (\beta), (\gamma)$ and
($\delta$) are defined as follows:
\begin{description}
\item{$(\alpha)$} if $P$ satisfies condition \textbf{U1}, then
operation $(\alpha)$ adds a new column made of $\ell(c_n)+1$ cells
on the right of $c_n$, see Figure~\ref{opA}.

\begin{figure}[h]
\begin{center}
%\centerline{\hbox{\includegraphics[width=0.9\textwidth]{operationA.eps}}}
\centerline{\hbox{\includegraphics[height=2.5cm]{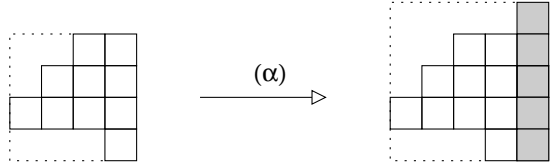}}}
\caption{Operation $(\alpha)$; the added column has been
highlighted}\label{opA}
\end{center}
\end{figure}
%***** INSERIRE FIGURA 7 di (1) *****

\item{$(\beta)$} it can be performed on each cell of $c_n$; so let
$d_i$ be the $i$th cell of $c_n$, from bottom to top, with $1 \leq
i \leq \ell(c_n)$. Operation $(\beta)$  adds a new row above the
row containing $d_i$ (of the same length), and add a new column on
the right of $c_n$ made of $i$ cells, see Figure~\ref{opB}.

\begin{figure}[h]
\begin{center}
%\centerline{\hbox{\includegraphics[width=0.9\textwidth]{operationB.eps}}}
\centerline{\hbox{\includegraphics[height=2.5cm]{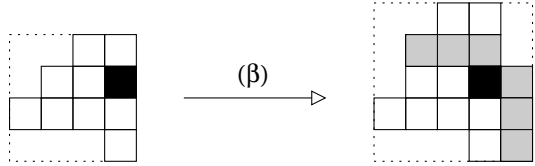}}}
\caption{Operation $(\beta)$; the cell $d_i$ is filled in black ,
the added column and row have been highlighted}\label{opB}
\end{center}
\end{figure}
%***** INSERIRE FIGURA 8 di (1) *****

\item{$(\gamma)$} it can be performed on each cell of $c_n$; so
let $d_i$ be the $i$th cell of $c_n$, from bottom to top, with $1
\leq i \leq \ell(c_n)$. Operation $(\gamma)$  adds a new row below
the row containing $d_i$ (of the same length), and add a new
column on the right of $c_n$ made of $n-i+1$ cells, see Figure
~\ref{opC}.
\begin{figure}[h]
\begin{center}
%\centerline{\hbox{\includegraphics[width=0.9\textwidth]{operationC.eps}}}
\centerline{\hbox{\includegraphics[height=2.5cm]{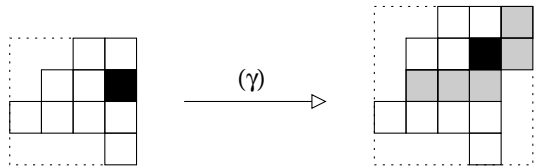}}}
\caption{Operation $(\gamma)$; the cell $d_i$ is filled in black ,
the added column and row have been highlighted}\label{opC}
\end{center}
\end{figure}
%***** INSERIRE FIGURA 9 di (1) *****

\item{$(\delta)$} if $P$ satisfies condition \textbf{U2}, then
operation $(\delta)$ adds a new column made of $\ell(c_n)+1$ cells
on the right of $c_n$, see Figure ~\ref{opD}.

\begin{figure}[h]
\begin{center}
%\centerline{\hbox{\includegraphics[width=0.9\textwidth]{operationD.eps}}}
\centerline{\hbox{\includegraphics[height=2.5cm]{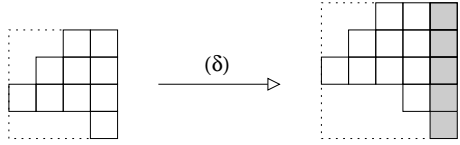}}}
\caption{Operation $(\delta)$; the added column has been
highlighted}\label{opD}
\end{center}
\end{figure}
%***** INSERIRE FIGURA 10 di (1) *****

Obviously in any case the obtained permutomino is a convex
permutomino of size $n+1$.

\end{description}

We refer to \cite{dfpr} for further details and proofs.

\section{The generating algorithm}\label{algo}
Our aim is to illustrate an exhaustive generating algorithm for
convex permutominoes, basing on the ECO construction recalled in
previous section. In the sequel we will refer only to convex
permutominoes, simply named ``permutomino".

First of all we define a subset of permutominoes of
size $n$, the so called {\em active permutominoes}, then we will
show the existence of a bijection between the active permutominoes
of size $n$ and the set of permutominoes  of size $n-1$. Finally,
we will define the generating tree of permutominoes of size $n$.

\subsection{Definition of active permutominoes}\label{sectactive}
A permutomino $P$ is {\em active} if the following conditions hold:
\begin{enumerate}
\item the leftmost column contains only one cell, ($\ell(c_1)=1$);
\item the leftmost reentrant point $\alpha$ has abscissa 2.
\end{enumerate}

\begin{figure}[h]
\begin{center}
%\centerline{\hbox{\includegraphics[width=0.9\textwidth]{active.eps}}}
\centerline{\hbox{\includegraphics[height=2.5cm]{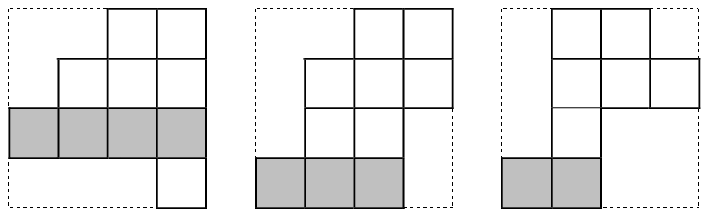}}}
\caption{Three active permutominoes of size $4$ }\label{active}
\end{center}
\end{figure}

\begin{figure}[h]
\begin{center}
%\centerline{\hbox{\includegraphics[width=0.9\textwidth]{not_active.eps}}}
\centerline{\hbox{\includegraphics[height=2.5cm]{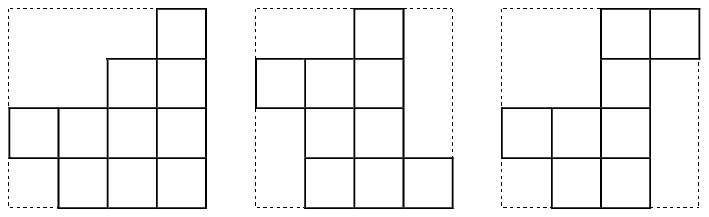}}}
\caption{Some not active permutominoes of size $4$}
\label{noactive}
\end{center}
\end{figure}

That is, the word encoding the boundary of an active
permutomino begins with $"NEN \ldots$ and ends with $\ldots WW"$.
In Figure~\ref{active} are depicted three active permutominoes of
size 4, while in Figure~\ref{noactive}   there are some
permutominoes that, yet having only one cell in the leftmost
column, are not active.

%**************** INSERIRE FIGURE 1 e 2 di (4) **********

\begin{figure}[h]
\begin{center}
%\centerline{\hbox{\includegraphics[width=0.9\textwidth]{not_admiss.eps}}}
\centerline{\hbox{\includegraphics[height=3.5cm]{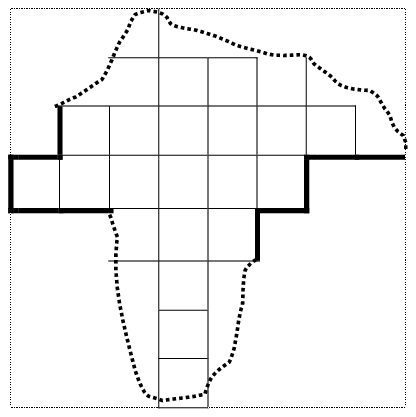}}}
\caption{A not admissible configuration } \label{noadmiss}
\end{center}
\end{figure}

Let $P$ be an active permutomino and let $\varrho$ be the row of
$P$ containing the only cell in the leftmost column (in
Figure~\ref{active} $\varrho$ is highlighted);
%$\varrho = r_i$, $1 \leq i < n$,
$\varrho$ can be the bottom row but it is never the top row. From
Proposition \ref{prop1}, $\varrho$ ends at the same abscissa of
either the above and the below row, if it exists; so a
configuration as that depicted in Figure~\ref{noadmiss} is not
admissible. Therefore, if in an active permutomino $P$ of size $n$
we remove $\varrho$, we obtain a convex permutomino $P'$ of size
$(n-1)$.

%**************** INSERIRE FIGURA  2bis di (4) **********

\begin{figure}[h]
\begin{center}
\centerline{\hbox{\includegraphics[width=0.9\textwidth]{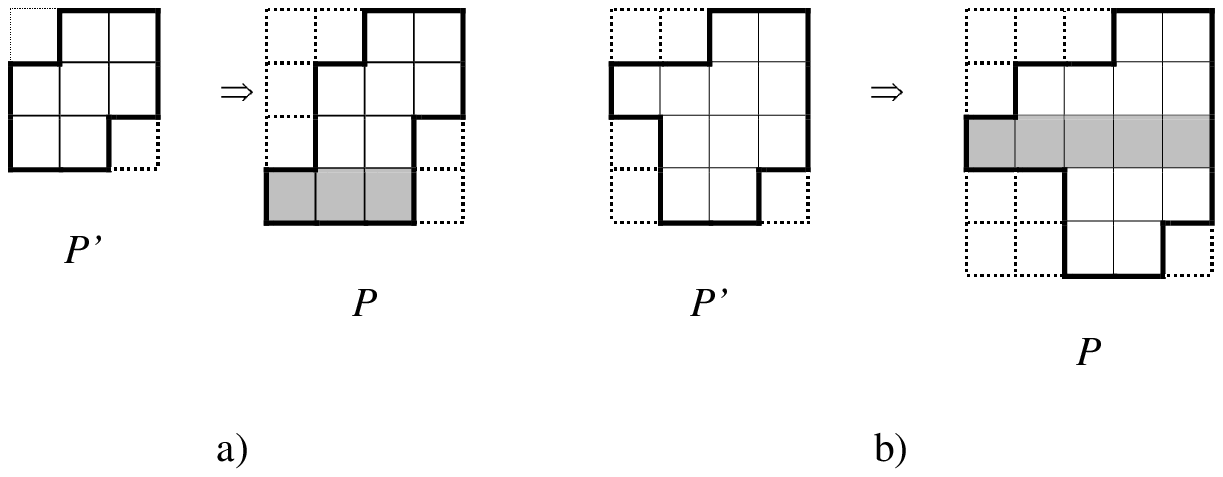}}}
\caption{A permutomino $P'$ of size $(n-1)$ and the corresponding
permutomino $P$ of size $n$} \label{from}
\end{center}
\end{figure}

Let $P'$ be a permutomino of size $(n-1)$ and let $\varrho'$ the
row in $P'$ containing the leftmost salient point with minimal
ordinate (point A in Figure~\ref{sample2}). If in $P'$ we add,
below $\varrho'$,  a  row $\varrho$ one cell longer on the left
and ending at the same abscissa of $\varrho'$, we obtain an active
permutomino $P$ of size $n$:
\begin{description}
\item{i.} if $\varrho'$ is the bottom row, the added row in $P'$
is one cell longer on the left and therefore the new permutomino
of size $n$ is convex and active (see Figure~\ref{from}~a));
\item{ii.} if $\varrho' = r_i$, with $1 < i <n$, the row $r_{i-1}$
ends at the same abscissa of $\varrho'$; so adding $\varrho$ we
obtain a convex and active permutomino of size $n$ (see Figure
~\ref{from}~b)).
\end{description}

This means that it exists a bijection
$\psi$ between the set ${\cal C}_{a,n}$ of active permutominoes and the set ${\cal
C}_{n-1}$:
$$
\psi \,: \, {\cal C}_{a,n} \, \to \, {\cal C}_{n-1}
$$
such that $\psi(P)$ is the permutomino of size $n-1$ obtained by
removing  $\varrho$ from the active permutomino $P$ of size $n$.

%**************** INSERIRE FIGURA 3 di (4) *******************

\subsection{The exhaustive generating algorithm}\label{gen}
The algorithm we propose for the exhaustive generation of convex
permutominoes of size $n$ is based on the bijection $\psi$,
defined in Section \ref{sectactive}, and on the ECO construction
of permutominoes recalled in Section \ref{eco}.

The generating process is described by an operator $\phi$ so
defined:\\

\begin{figure}[h]
\begin{center}
%\centerline{\hbox{\includegraphics[width=0.9\textwidth]{firstperm.eps}}}
\centerline{\hbox{\includegraphics[height=2.5cm]{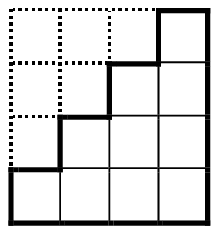}}}
\caption{$P_n^r$ of size 4} \label{first}
\end{center}
\end{figure}

\textbf{$\phi$ Operator:}
\begin{enumerate}
\item The first permutomino of the generating process is
$P_{n}^r$, that is the permutomino of size $n$ associated with the
pair of permutations $(\pi_1,\pi_2)$ of $[n+1]$:
$$
\pi_1 = (1,2,3, \ldots, n, n+1)\qquad
\pi_2 = (2,3,4, \ldots, n+1,1)
$$
In Figure~\ref{first} is depicted $P_n^r$  of size 4.

%**************** INSERIRE FIGURA 3 bis di (4) *******************

\item $P_{n}^r$
is an active permutomino; let $\bar{P}^r$ the permutomino of size
$(n-1)$ such that
 ${\bar P}^r= \psi(P_{n}^r)$ ({\em $\bar{P}^r$  is obtained by removing the bottom row of $P_{n}^r$ }).
\item Apply operations $(\beta)$, $(\gamma)$ and $(\delta)$ of ECO
construction to $\bar{P}^r$. Every application generates a new
convex permutomino of size $n$. \item For each new active
permutomino $Q$ repeat the following actions until active
permutominoes are generated:
\begin{description}
\item{4.1} remove the $\varrho$ row from $Q$ obtaining $\bar{Q} = \psi(Q)$;
\item{4.2} apply all the possible operations of the ECO construction to $\bar{Q}$. Every
application generates a new convex permutomino of size $n$.
\end{description}
\end{enumerate}

Our strategy can be represented using a rooted tree, say {\em
${\cal C}_n$-tree}, so defined:\\
%\textbf{${\cal C}_n$-tree}:
\begin{enumerate}
\item the root is $P_{n}^r$ and it is at level 0; \item if $Q \in
{\cal C}_n${\em-tree} is an active permutomino at level $k \geq
0$, then $\phi(Q) = \vartheta(\psi(Q))$ ($\vartheta$ is the
operator defined in the ECO construction) and every $P \in
\phi(Q)$ is a son of $Q$ and it is at level $(k+1)$. For the sake
of simplicity, we say that $P \in \phi^{k+1}(P_n^r)$.
\end{enumerate}

\begin{figure}[p]
\begin{center}
%\centerline{\hbox{\includegraphics[width=0.9\textwidth]{firstperm.eps}}}
\centerline{\hbox{\includegraphics[height=15cm]{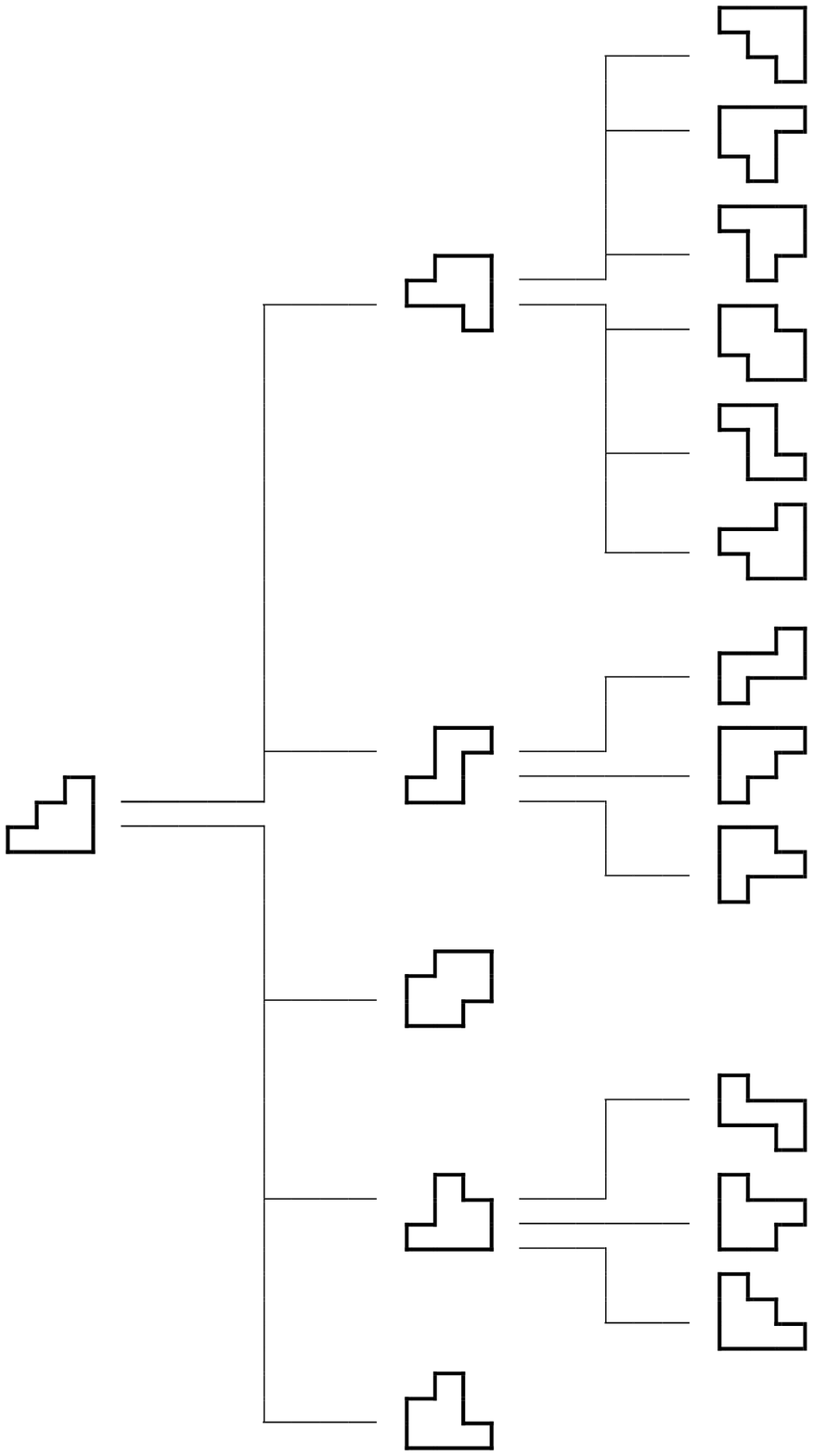}}}
\caption{${\cal C}_3${\em-tree}} \label{tree3}
\end{center}
\end{figure}

In Figure~\ref{tree3} ${\cal C}_3${\em-tree} is illustrated.

%**************** INSERIRE FIGURA 5 di (4) *******************

\begin{proposition}
${\cal C}_n${\em-tree} contains all and only the
convex permutominoes of size $n$, i.e ${\cal C}_n${\em-tree} $=
{\cal C}_n$.
\end{proposition}
\noindent \textbf{Proof}. \textbf{[Only]} The permutominoes in ${\cal
C}_n${\em-tree} are obtained from permutominoes of size $(n-1)$ by
applying the ECO construction; so, as proved in \cite{dfpr}, we
generate permutominoes of size $n$. Therefore, since the
permutominoes of size $(n-1)$ are each other different, the ones
of size $n$ are different too.

\noindent  \textbf{[All]} We must proof that for each permutomino
$Q \in {\cal C}_n$ there exists a path from $P_n^r$ to $Q$, that
is  it exists a finite sequence $P_0, P_1, \ldots, P_k$ with  $k
\in \mathbb{N}$ and $P_k=Q$ such that:
\begin{itemize}
    \item $P_0=P^n_{r}$;
    \item $P_{i+1} \in \phi(P_i)$, \hspace{0.5cm}$0\leq$ i $\leq$
    $k-1$.
\end{itemize}

\begin{figure}[h]
\begin{center}
%\centerline{\hbox{\includegraphics[width=0.9\textwidth]{firstperm.eps}}}
\centerline{\hbox{\includegraphics[height=3cm]{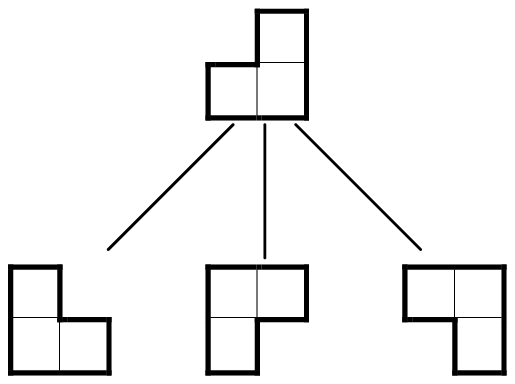}}}
\caption{${\cal C}_{2}${\em-tree}} \label{tree2}
\end{center}
\end{figure}

In other words, we must prove that there exists $k \geq 1$ such
that $Q \in \phi^k(P^r_n)$. We know that the active permutominoes
of size $n$ are as many as the permutominoes of size $(n-1)$. So
it is sufficient to proof the following:
\begin{proposition}\label{prop2}
All the active permutominoes of size $n$ are generated.
\end{proposition}
\noindent  \textbf{Proof}. By induction on the size $n$.

\noindent  {\em Base.} For $n = 1$ there is only the permutomino
containing one cell; if $n = 2$ the unique active permutomino is
$P^r_2$, (see Figure~\ref{tree2}). So for $n \leq 2$ Proposition
\ref{prop2} yields.

\noindent  {\em Inductive hypothesis.} Let us assume that all
the permutominoes of size $(n-1)$ are generated and let ${\cal C}_
 {n-1}${\em-tree} be the associated tree. Then, starting from the root
$P_{n-1}^r$ it is possible to reach any permutomino of size
$(n-1)$. So, for each permutomino $P_{n-1}\in {\cal C}_{n-1}$ there
exists a $k$ such that:
$$
P_{n-1} \in \phi^k(P_{n-1}^r)
$$
\noindent  {\em Inductive step.} Let $\bar{P}_n$ be an active
permutomino of size $n$ and let $Q_{n-1} = \psi(\bar{P}_n)$
($Q_{n-1}$ is obtained from $\bar{P}_n$ removing  $\varrho$). So
there exists $\bar{k}$ such that
$$
Q_{n-1} \in \phi^{\bar k}(P_{n-1}^r).
$$
But
$$
P_{n-1}^r = \psi(P_{n}^r)
$$
so it follows that:
$$
Q_{n-1} \in \phi^{\bar k}(P_{n}^r)
$$
that is, there is a path from the root $P_{n}^r$ to any
permutomino of size $(n-1)$. Therefore,  each permutomino of size
$n$ is reachable from $P_{n}^r$.\hspace{2.0cm}$\square$\\
%\begin{flushright}
%$\Box$
%\end{flushright}

%**************** INSERIRE FIGURA 4 di (4) *******************

\subsection{Algorithm cost analysis}

First of all, we will prove that the height of ${\cal
C}_{n}${\em-tree} is $n$. The proof is based on the following
propositions.
\begin{proposition}\label{prop3}
Using the ECO construction of Section \ref{eco}, an active
permutomino of size $n$ is generated by one and only one active
permutomino of size $(n-1)$.
\end{proposition}
\noindent  \textbf{Proof}. It follows straightforward from the ECO
construction which never adds a column on the left containing one cell to the permutoninoes.
\hspace{1.0cm}$\square$\\
%\begin{flushright}
%$\Box$
%\end{flushright}

\begin{proposition}\label{prop4}
Given an active permutomino $Q$, the longest path starting from
$Q$ has length $j$ if its $\alpha$ points lie in $(2,h),
(3,h+1), \ldots, (j+1,h+j-1)$, $h$ being the ordinate of the lefmost $\alpha$ point.
\end{proposition}
\noindent  \textbf{Proof}. The permutomino $\psi(Q)$ of size
$(n-1)$ to which the ECO construction is applied, is obtained
removing from $Q$ the row $\varrho$, so the leftmost $\alpha$
point of $Q$ is removed in $\psi(Q)$. Therefore, $\psi(Q)$ will be
active, and then, from Proposition \ref{prop3}, it can generate
new active permutominoes, only if it has an $\alpha$ point at
abscissa 2; thus, $Q$ must have an $\alpha$ point at abscissa 3.
In the same way, a permutomino generated from $\psi(Q)$ will be
active only if $\psi(Q)$ has an $\alpha$ point at abscissa 3, that
is if $Q$ has an $\alpha$ point at abscissa 4, and so on up to
$(j+1)$.\hspace{9.0cm}$\square$\\
%\begin{flushright}
%$\Box$
%\end{flushright}

The permutomino $P_{n}^r$ has $(n-1)$ consecutive $\alpha$ points
in $(2,2), \ldots, (n,n)$, so, from Proposition
\ref{prop4}, the longest path starting from  $P_{n}^r$ has length $(n-1)$.
Therefore, the height of ${\cal C}_{n}${\em-tree} is $n$.

From the generating alghorithm of Section \ref{gen} it follows
that the generation of a permutomino $P$ at level $k$ in the
${\cal C}_ {n}${\em-tree} depends only on a permutomino at level
$(k-1)$, that is, ${\cal C}_ {n}${\em-tree} is generated level by
level. Therefore, since the generation of a permutomino of size
$n$ from one of size $(n-1)$ has a constant cost, we may conclude
that the cost of the exhaustive generating algorithm is
proportional to the number of permutominoes of size $n$.

%\begin{thebibliography}{13}

\end{document}